\documentclass[12pt]{amsart}
\usepackage{url, amssymb,amsmath,amstext, amsthm, amsfonts, amssymb}
\usepackage{mathrsfs}
\usepackage{color}

\usepackage{hyperref}
\usepackage{cite}
\usepackage[all]{xypic}

\usepackage{graphicx}
\newtheorem{theorem}{Theorem}[section]

\newtheorem{prop0}[theorem]{Proposition}
\newtheorem{thm0}[theorem]{Theorem}
\newtheorem{conj0}[theorem]{Conjecture}
\newtheorem{lemma0}[theorem]{Lemma}

\newtheorem{computation}[theorem]{Computation}
\theoremstyle{remark}
\newtheorem{remark}[theorem]{Remark}

\theoremstyle{definition}

\def\rig#1{\smash{ \mathop{\longrightarrow}
    \limits^{#1}}}

\def\O{{\mathcal O}}

\def\o{{\otimes}}

\newcommand{\Hom}{\operatorname{Hom}}

\newcommand{\Osh}{{\mathcal O}}

\newcommand{\rank}{\operatorname{rank}}

\newcommand{\PP}{\mathbb{P}}
\newcommand{\CC}{\mathbb{C}}

\newcommand{\ZZ}{\mathbb{Z}}

\newcommand{\sL}{\mathcal{L}}
\newcommand{\sM}{\mathcal{M}}
\newcommand{\sH}{\mathcal{H}}
\newcommand{\sA}{\mathcal{A}}
\newcommand{\Sym}{\operatorname{Sym}}

\newcommand{\Seg}{\operatorname{Seg}}

\def\bw#1{{\textstyle\bigwedge^{\hspace{-.2em}#1}}}
\def\phi{\varphi}

\date{\today}

\begin{document}
\title[Homotopy techniques \& tensors]{Homotopy techniques for tensor decomposition \\ and perfect identifiability}

\author[J.D. Hauenstein]{Jonathan D. Hauenstein}
\address{Department of Applied and Computational Mathematics and Statistics, University of Notre Dame, Notre Dame, IN, USA}
\email{hauenstein@nd.edu}
\urladdr{\url{http://www.nd.edu/~jhauenst}}

\author[L. Oeding]{Luke Oeding}
\address{Department of Mathematics and Statistics, Auburn University,
Auburn, AL, USA}
\email{oeding@auburn.edu}
\urladdr{\url{http://www.auburn.edu/~lao0004}}

\author[G. Ottaviani]{Giorgio Ottaviani}
\address{Dipartimento di Matematica e Informatica ``U. Dini'' \\
Universit\`a degli Studi di Firenze \\
Firenze, Italy}
\email{ottavian@math.unifi.it}
\urladdr{\url{http://web.math.unifi.it/users/ottaviani}}

\author[A.J. Sommese]{Andrew J. Sommese}
\address{Department of Applied and Computational Mathematics and Statistics, University of Notre Dame, Notre Dame, IN, USA}
\email{sommese@nd.edu}
\urladdr{\url{http://www.nd.edu/\~sommese}}

\maketitle

\begin{abstract}
Let $T$ be a general complex tensor of format $(n_1,...,n_d)$. When the fraction $\prod_in_i/[1+\sum_i(n_i-1)]$ is an integer,
and a natural inequality (called balancedness) is satisfied, it is expected that $T$ has finitely many minimal decomposition as a sum
of decomposable tensors.
We show how homotopy techniques allow us to find all the
decompositions of~$T$, starting from a given one.
Computationally, this gives a guess
regarding the total number of such decompositions.
This guess matches exactly with all cases previously known, and predicts several unknown cases.
Some surprising experiments yielded two new cases of
generic identifiability: formats $(3,4,5)$ and $(2,2,2,3)$
which have a unique decomposition as the sum of $6$ and $4$ decomposable tensors, respectively.  We conjecture that these two cases
together with the classically known matrix pencils are the only cases where generic identifiability holds, i.e., the only {\em identifiable} cases.
Building on the computational
experiments, we use algebraic geometry
to prove these two new cases are~indeed~generically~identifiable.
\end {abstract}

\section{Introduction}
Tensor decomposition is an active field of research, with many applications
(see, e.g.,~\cite{Land12} for a broad overview).
A tensor $T$ of format $(n_1,\ldots, n_d)$ is an element of the
tensor space $\CC^{n_1}\otimes\cdots\otimes\CC^{n_d}$.
The {\it rank} of $T$ is the minimum $r$ such that
\begin{equation}\label{tensordec:eq}
T=\sum_{i=1}^r v_1^i\otimes\cdots\otimes v_d^i
\end{equation}
where $v_j^i\in\CC^{n_j}$.
This reduces to the usual matrix rank when $d = 2$.

The space $\CC^{n_1}\otimes\cdots\otimes\CC^{n_d}$
contains a dense subset where the rank is constant.
This constant is called the {\it generic rank}
for tensors of format $(n_1,\dots,n_d)$.
By a simple dimensional count,
the generic rank for tensors of format $(n_1,...,n_d)$ is at least
\begin{equation}\label{eq:R}
R(n_1,\dots,n_d) := \frac{\prod_{i=1}^dn_i}{1+\sum_{i=1}^d(n_i-1)}
= \frac{\prod_{i=1}^dn_i}{1-d + \sum_{i=1}^d n_i}.
\end{equation}
The value $\lceil R(n_1,\dots,n_d)\rceil$ is called
the {\it expected generic rank} for $(n_1,\dots,n_d)$.

A necessary condition for a general tensor $T$ of format $(n_1,\dots,n_d)$ to have only finitely many decompositions (\ref{tensordec:eq}) is that the number $R(n_1,\dots,n_d)$ is actually an integer.
Such formats are called {\it perfect} \cite{Stra, BCS}.
Moreover, if a general tensor is known to have finitely
many decompositions (\ref{tensordec:eq}),
then the generic rank is equal to
the expected generic rank $R(n_1,\dots,n_d)$.

 {A tensor format is said to be \emph{generically identifiable} if a generic tensor of that format has a unique decomposition (up to global rescaling and permutation of the factors). A tensor format is said to satisfy \emph{perfect identifiability} if the format is perfect and generically identifiable.}
The main goal of this paper is to study the number
of decompositions of perfect formats $(n_1,\dots,n_d)$
when the generic rank is indeed equal to the generic expected rank.

The main tool for inquiry is numerical algebraic geometry,
a collection of algorithms to numerically compute and manipulate
solutions sets of polynomial systems.
Numerical algebraic geometry, named in \cite{SW95}, grew out of
numerical continuation methods for finding all isolated solutions
of polynomial systems.
For a development and history of the area, see the monographs \cite{SWbook,Bertini} and the survey \cite{WS11}.
The monograph \cite{Bertini} develops the subject using
the software package {\tt Bertini} \cite{BertiniSoftware},
which is used to perform the computations in this article.
For understanding the relation between numerical approaches
and the more classical symbolic approaches to
computational algebraic geometry, see \cite{BDHPPSSW}.

Numerical algebraic geometry has proven useful in many other
applications.  A small subset of such applications include
computing the initial cases for equations of an infinite family
of Segre-Grassmann hypersurfaces in \cite{DHO};
numerically decomposing a variety in \cite{BatesOeding}
which was a crucial computation leading to a set-theoretic solution of
the so-called {\it salmon problem} \cite{AllmanPrize}
improving upon a previous result of Friedland \cite{Fri13},
and inspiring a later result of Friedland and Gross \cite{FG12};
solving Alt's problem \cite{Alt} which counts
the number of distinct four-bar linkages whose coupler curve
interpolates nine general points in the plane,
namely 1442 \cite{NinePoint};
finding the maximal likelihood degree
for many cases of matrices with rank constraints \cite{MLDegree}
and observing duality which was proven in \cite{DR};
a range of results in physics such as \cite{HNS,HNS2,MDHS,MKH};
and numerically solving systems of nonlinear differential equations \cite{HHHS,HHHLSZ,HHHLSZ2,HHHLSZ3,HHHMS,HHSSXZ,HHS,HHS2}.

We consider the equation (\ref{tensordec:eq})
where $r$ is the generic rank and the $v_j^i$'s are unknowns.
Starting from one decomposition for $T$, we can move $T(s)$ along a
loop, for $0\le s\le 1$, such that $T(0)=T(1)=T$.
This consequently defines corresponding vectors $v_j^i(s)$ which satisfy
$$T(s)=\sum_{i=1}^r v_1^i(s)\otimes\cdots\otimes v_d^i(s).$$

{The decompositions of $T$ at the end values, $s = 0$ and $s = 1$,
may be different.}  Since this process is computationally cheap,
it can be repeated with random loops a considerable number of times
and one can record all of the distinct decompositions found.
Moreover, in the perfect case,
where decompositions correspond to solutions
to system of polynomial equations with the same number
of variables, i.e., a square system, one can
use $\alpha$-theory via {\tt alphaCertified}\cite{alphaCertified,alphaCertifiedSoftware}
to prove lower bounds on the total number of decompositions.
{Theory guarantees that all decompositions can be found
using finitely many loops while}
experience shows that
{all decompositions of $T$ can be found after
a certain number of attempts using random loops.}
When the number of decompositions is small,
this process stabilizes quickly to the total number
of decompositions.
We describe using this process on some previously known
cases and predict several unknown cases.
In particular, the values reported in Section~\ref{sec:Homotopy}
are provable lower bounds that {we expect are sharp}.

To put these results in perspective, we recall that finding equations that detect tensors of small rank is a difficult subject.
Recent progress is described in \cite{AJRS}, which gives a semi-algebraic description of tensors of format $(n,n,n)$ of rank $n$ and
multilinear rank $(n,n,n)$.
In addition, several recent algorithms and techniques are
available to find best rank-one approximations \cite{NieLi}
or even to decompose a tensor of small
rank \cite{Tensorlab, BCMT, BMV, OeOtt}.
However, the problem generally becomes more difficult as the
rank increases so that decomposing a tensor which has the
generic rank is often the hardest case.

The formats $(3,4,5)$ and $(2,2,2,3)$ were exceptional in our
series of experiments since our technique showed that they
have a unique decomposition (up to reordering).
Indeed, an adaptation of the approach
developed in \cite{OeOtt} allowed to us to confirm our computations.

\begin{thm0}\label{345:thm}
A general tensor of format $(3,4,5)$ has a unique decomposition (\ref{tensordec:eq})
as a sum of $6$ decomposable summands.
\end{thm0}

\begin{thm0}\label{thm2223}
A general tensor of format $(2,2,2,3)$
has a unique decomposition (\ref{tensordec:eq})
as a sum of $4$ decomposable summands.
\end{thm0}

These theorems are proved in Section~\ref{sec:Eigen}.
The proofs provide algorithms for computing the unique decomposition,
which we have implemented in {\tt Macaulay2} \cite{M2}.
Based on the evidence described throughout,
we formulate the following conjecture.

\begin{conj0}\label{conj:Gen}
The only perfect formats $(n_1,\dots,n_d)$, i.e., $R(n_1,\dots,n_d)$
in \eqref{eq:R} is an integer, where a general tensor
has a unique decomposition \eqref{tensordec:eq} are:
\begin{enumerate}
\item $(2,k,k)$ for some $k$ --- matrix pencils, known classically by Kronecker normal form,
\item $(3,4,5)$, and
\item $(2,2,2,3)$.
\end{enumerate}
\end{conj0}

We would like to contrast the tensor case to the symmetric tensor case,
where the exceptional cases were known since the $19^{\rm th}$ century, 
as well as the following recent result.

\begin{thm0}[Galuppi--Mella~\cite{GM}]\label{conj:Sym}
The only perfect formats $(n,d)$, i.e., $n^{-1}\cdot{{n+d-1}\choose d}$ is an integer, where a general tensor in $\Sym^d\CC^n$
has a unique decomposition are:
\begin{enumerate}
\item $(2,2k+1)$ for some $k$ --- odd degree binary forms, known to Sylvester,
\item $(3,5)$ --- Quintic Plane Curves (Hilbert, Richmond, Palatini), and
\item $({4},3)$ --- Cubic Surfaces (Sylvester Pentahedral Theorem).
\end{enumerate}
\end{thm0}
Theorem  \ref{conj:Sym} was still stated  as a conjecture 
(following \cite{Mella06} and \cite{Mella09}) in the first preprint of the present article.
See \cite{Mella06} and \cite{RS} for classical references. In \cite[\S~4.4]{OeOtt},
two of the authors showed that the Koszul flattening method predicts
exactly the cases listed in Theorem \ref{conj:Sym} and no others. 

\section{Some known results on the number of tensor decompositions}\label{sec:history}
\subsection{General tensors}

The following summarizes some known results about tensors
of format $(n_1,\ldots,n_d)$.
For any values of $r$ smaller than the generic rank,
which was defined in the introduction,
the (Zariski) closure of the set of tensors of
rank $r$ is an irreducible algebraic variety.
This variety is identified with the cone over
the $r^{\rm th}$ secant variety to the Segre variety $\PP(\CC^{n_1})\times\cdots\times\PP(\CC^{n_d})$
of decomposable tensors, e.g., see \cite{Land12, CGO}.
In particular, it is meaningful to speak about a general tensor
of rank $r$.

Throughout this section, we consider cases where $d\geq3$
and, without loss of generality, assume that
$2\le n_1\le n_2\le\ldots \le n_d$.
First, we review the known results on the so-called unbalanced formats.

\begin{thm0}\label{known_general:thm}
For formats $(n_1,\dots,n_d)$, suppose that
$n_d\ge\prod_{i=1}^{d-1}n_i-\sum_{i=1}^{d-1}(n_i-1)$.
\begin{enumerate}
\item\label{gen:Item1} The generic rank is $\min\left(n_d,~\prod_{i=1}^{d-1}n_i\right)$.\smallskip

\item\label{gen:Item2} A general tensor of rank $r$ has a unique decomposition if
$r<\prod_{i=1}^{d-1}n_i-\sum_{i=1}^{d-1}(n_i-1)$.\medskip

\item\label{gen:Item3} A general tensor of rank $r=\prod_{i=1}^{d-1}n_i-\sum_{i=1}^{d-1}(n_i-1)$ has exactly $D\choose r$ different decompositions
 where
$$D=\frac{\left(\sum_{i=1}^{d-1}(n_i-1)\right)!}{(n_1-1)!\cdots(n_{d-1}-1)!}.$$
This value of $r$ coincides with the generic rank
in the perfect case: when $r = n_d$.\medskip

\item\label{gen:Item4} If $n_d > \prod_{i=1}^{d-1}n_i-\sum_{i=1}^{d-1}(n_i-1)$,
a general tensor of rank $r>\prod_{i=1}^{d-1}n_i-\sum_{i=1}^{d-1}(n_i-1)$,
e.g., a general tensor of format $(n_1,\dots,n_d)$, has infinitely many decompositions.
\end{enumerate}
\end{thm0}
\begin{proof}
When $n_d>\prod_{i=1}^{d-1}n_i-\sum_{i=1}^{d-1}(n_i-1)$,
Item~\ref{gen:Item1} follows from \cite[Thm.~2.1(1-2)]{CGG02} (see also \cite[Prop.~8.2]{BCO}).
In the perfect case, i.e.,
$n_d=\prod_{i=1}^{d-1}n_i-\sum_{i=1}^{d-1}(n_i-1)$,
Item~\ref{gen:Item1} follows from \cite[Prop.~8.3]{BCO}.
Items~\ref{gen:Item2} and~\ref{gen:Item3} follow from
\cite[Prop.~8.3,~Cor.~8.4]{BCO}.
When $n_d-1>\prod_{i=1}^{d-1}n_i-\sum_{i=1}^{d-1}(n_i-1)$,
Item~\ref{gen:Item4} follows from \cite[Lemma 4.1]{AOP}.
If $n_d-1=\prod_{i=1}^{d-1}n_i-\sum_{i=1}^{d-1}(n_i-1)$,
then $\prod_{i=1}^{d-1}n_i=\sum_{i=1}^{d}(n_i-1)$.
Hence, $1+\sum_{i=1}^{d}(n_i-1)$ cannot divide
$\prod_{i=1}^{d}n_i$ and so the format cannot be perfect.
\end{proof}

The case $(2,n,n)$, corresponding to pencils of square matrices, is the only case for which the binomial coefficient
${D\choose r}$ in Theorem~\ref{known_general:thm}(\ref{gen:Item3}) is equal
to $1$.  The unique decomposition is a consequence of the canonical form
for these pencils, found by Weierstrass~and~Kronecker~\cite{BCS}.

For convenience, {Table~\ref{table2.1} lists} some perfect cases coming from Theorem \ref{known_general:thm}(\ref{gen:Item3}), namely
when $n_d=\prod_{i=1}^{d-1}n_i-\sum_{i=1}^{d-1}(n_i-1)$
with generic rank $r = n_d$.

\begin{table}
\[\begin{array}{|r|r|r|}
\hline
(n_1,\ldots,n_d)&\textrm{gen. rank}&\textrm{\# of decomp. of general tensor}\\
\hline
(2, n,n)& n& 1\,\,\\
\hline
(3, 3, 5)& 5& 6\,\,\\
\hline
(3, 4, 7)& 7& 120\,\,\\
\hline
(3, 5, 9)& 9& 5005\,\,\\
\hline
(3, 6, 11)& 11& 352716\,\,\\
\hline
(4, 4, 10)& 10& 184756\,\,\\
\hline
(2, 2, 2, 5)& 5& 6\,\,\\
\hline
(2, 2, 3, 8)& 8& 495\,\,\\
\hline\end{array}\]
\caption{Generic ranks and numbers of decompositions for perfect formats (Thm.~\ref{known_general:thm}).}\label{table2.1}
\end{table}
\smallskip

After Theorem \ref{known_general:thm}, the only open cases
are when the balancedness condition is satisfied:
\begin{equation}\label{balanc:eq}n_d<\prod_{i=1}^{d-1}n_i-\sum_{i=1}^{d-1}(n_i-1).\end{equation}

A seminal identifiability result for general tensors up to a
certain rank is \cite[Cor.~3.7]{Stra}.
In \cite{COV2014}, based on weak defectiveness introduced
in \cite{CC}, there are techniques to check the number of
decompositions of a general tensor of rank $r$,
generalizing Kruskal's result~\cite{Kru}.

De Lathauwer's condition (mentioned after equation (1.7) in \cite{DeLathauwer}) can guarantee uniqueness up to rank
\begin{equation}\label{eq:LDL}
 \sum_{i=1}^{d-1} n_i - d+1
\end{equation}
In the $3\times 4\times 5$ the bound \eqref{eq:LDL} is  $5$.
In the $2\times 2\times 2\times 3$ the bound \eqref{eq:LDL} is  $3$.  
Indeed De Lathauwer considers only ``tall'' tensors when the rank is no larger than the dimension of any mode, so his methods don't apply to the perfect case. However if you ignore this assumption (rank $R\leq n_{3}$) in  \cite[Theorem~2.5]{DeLathauwer}, the second inequality 
\[
R(R-1) \leq \frac{n_{1}(n_{1}-1)n_{2}(n_{2}-1)}{2}
\]
becomes $30 \leq 36 $, which is valid in the $3\times 4\times 5$ case. So the previous bound just misses our results.

For all formats such that $\prod_{i=1}^dn_i\le 15,000$
which satisfy the inequality (\ref{balanc:eq}),
a general tensor of rank $r$ which is strictly
smaller then the generic rank has a unique decomposition
except for a list of well understood exceptions, e.g., see \cite[Thm.~1.1]{COV2014}.
These results support the belief that, other than
some exceptions, a general tensor of rank $r$
smaller then the generic rank has a unique decomposition.
When $r$ is the generic rank,
since the techniques in \cite{COV2014} cannot be applied,
we apply numerical algebraic geometry
to such cases in Section~\ref{sec:Homotopy}.

\subsection{The symmetric case}

The following summarizes results about symmetric tensors to contrast
with the general case.
Recall that symmetric tensors of format $(n,d)$
are tensors $T\in\Sym^d\CC^{n}$, which can be identified with homogeneous polynomials of degree $d$ in $n$ variables.
The (symmetric) rank of $T$ is the minimum $r$ such that
there is an expression
$$T=\sum_{i=1}^r v^i\otimes\cdots\otimes v^i$$
with $v^i\in\CC^{n}$. If $T$ is identified with a polynomial, then each summand $v^i\otimes\cdots\otimes v^i$
is the $d$-power of a linear form.
By a naive dimension count, a general tensor
in $\Sym^d\CC^{n}$ has rank at least $n^{-1}\cdot{{n+d-1}\choose d}$.
When this fraction is an integer, the symmetric format $(d, n)$ is called {\it perfect}. As in the general case, perfectness is a necessary condition
in order for a general tensor in $\Sym^d\CC^{n}$
to have only finitely many decompositions.

The following is the basic result about decomposition of symmetric
tensors that we state for perfect formats.

\begin{thm0}[Alexander-Hirschowitz\cite{AH}]\label{AH}
Let $d\ge 3$ and assume \mbox{$r=n^{-1}\cdot {{n+d-1}\choose {d}}\in\ZZ$}.
A~general tensor in $\Sym^d\CC^{n}$ has
finitely many decompositions of rank $r$ except in the following cases:
\mbox{$(n,d) = (3,4),~(5,3),~\mbox{or~}(5,4)$}.
In these three exceptional cases, a general tensor has no
decomposition of rank $r$, but infinitely many decompositions of rank $r+1$.
\end{thm0}

When $n=2$, note that $\frac{d+1}{2}\in\ZZ$ exactly when $d$ is odd.
In these cases, Sylvester proved that there is a unique
decomposition with $\frac{d+1}{2}$ summands \cite{RS}.

Theorem~\ref{conj:Sym}(\cite{GM}) lists
two other cases when a general tensor
in $\Sym^d\CC^{n}$ has a unique decomposition,
namely $\Sym^3\CC^{4}$ {where $(n,d) = (4,3)$}
and $\Sym^5\CC^{3}$ {where $(n,d) = (3,5)$} \cite{RS}.

When $n=3$, note that $\frac{1}{3}\cdot{{d+2}\choose {2}}\in\ZZ$
is an integer exactly when $d=1$ or $d=2$ modulo~$3$.
{Table~\ref{table2.2} records all the cases that can be found in \cite{RS} concerning $\Sym^d\CC^{3}$.}
\begin{table}
\[\begin{array}{|c|c|c|}
\hline
d&\textrm{gen. rank}&\textrm{\# of decomp. of general tensor}\\
\hline
4& 6& \infty\\
\hline
5& 7& 1\\
\hline
7& 12& 5\\
\hline
8& 15& 16\\
\hline
\end{array}\]
\caption{Generic ranks and numbers of decompositions for general tensors in $\Sym^d\CC^{3}$ \cite{RS}.}\label{table2.2}
\end{table}
The clever syzygy technique used in \cite{RS} seems not to
extend to higher values of $d$.

\begin{remark}\label{rem:tangproj}
Let $d=1$ or $d=2$ modulo~$3$. By \cite[Thm.~4.2(vi)]{CR}, the number
of decompositions of a general symmetric tensor in $\Sym^d\CC^{3}$  is bounded below by the degree
of the tangential projection from  $r-1$ points, where $r=\frac{(d+2)(d+1)}{6}$ is the generic rank.
This latter degree is computed as the residual intersection
of two plane curves of degree $d$ having $r-1$ double points, which is
$d^2-4(r-1)=d^2-4\frac{(d-1)(d+4)}{6}=\frac{(d-2)(d-4)}{3}$.

An analysis of the degeneration performed in \cite{CR} suggests that actually
the number
of decompositions of general symmetric tensor in $\Sym^d\CC^{3}$
  should be divisible by
$\frac{(d-2)(d-4)}{3}$.
 This guess agrees, for $d\le 8$, with the above table from \cite{RS} and
the results for $d\le 11$~in~\mbox{\S \ref{subsec:symmetric}}.
\end{remark}

A general symmetric tensor of rank $r$ which is strictly
smaller then the generic rank has a unique decomposition
except for a list of well understood exceptions, see \cite{Bal, Mella06, COV2}.

\section{Homotopy techniques for tensor decomposition}\label{sec:Homotopy}

In this section, we first describe the monodromy-based approach
we use to determine the number of decompositions for a general tensor.
The software {\tt Bertini}~\cite{BertiniSoftware,Bertini}
is then used in the subsequent subsections to compute
decompositions for various formats.  In the perfect cases
under consideration, the number of decompositions can be
certifiably lower bounded via {\tt alphaCertified}~\cite{alphaCertified,alphaCertifiedSoftware}.
In particular, for the two cases $(3, 4, 5)$ and $(2,2,2,3)$
which are discovered here to have a unique decomposition,
we provide theoretical proofs in Section~\ref{sec:Eigen}.
{Our computational methods include probabilistic reductions (e.g. cutting by random hyperplanes and choosing random points) and numerical computations which are always subject to round-off errors in any finite computation (e.g. numerical path tracking). Even though these methods are now completely standard, have been carefully and repeatedly tested, and yield completely reproducible results, they are technically only true with high probability, or up to the numerical precision of the computers we use. Recently results whose proofs partially rely on such methods have been denoted Theorem$^{\star}$ (see \cite{OedingSam}). In this article we say that these results hold \emph{with high confidence}.
}
\subsection{Decomposition via monodromy loops}

In numerical algebraic geometry, monodromy loops
have been used to decompose solution sets into irreducible
components~\cite{monodromy}.  Here, we describe the use of
monodromy loops for computing additional decompositions
of a general tensor.  For demonstration purposes, suppose
that a general tensor of format $(n_1,\dots,n_d)$ has
rank $r$ and finitely many decompositions.

The approach starts with a general tensor $T$
of format $(n_1,\dots,n_d)$ with a known decomposition
\eqref{tensordec:eq} with $v_j^i\in\CC^{n_j}$ for $i = 1,\dots,r$.
In practice, one randomly selects the $v_j^i$ first
and then computes the corresponding $T$ defined
by \eqref{tensordec:eq}.
To remove the trivial degrees of freedom, we assume
that $(v_j^i)_1 = 1$ for $i = 1,\dots,r$ and $j = 1,\dots,d-1$.
That is, we have a solution of
$$F_T(v_1^1,\dots,v_d^r) = \left[\begin{array}{cl}
T - \sum_{i=1}^r v_1^i\otimes\cdots\otimes v_d^i & \\
(v_j^i)_1 - 1, & i=1,\dots,r, j=1,\dots,d-1 \end{array}\right] = 0.$$
The system $F_T$ consists of $\prod_{j=1}^d n_j + r(d-1)$
polynomials in $r\cdot\sum_{j=1}^d n_j$ variables.
Since $r = R(n_1,\dots,n_d)$ in \eqref{eq:R},
the number of polynomials is equal to the number of variables
meaning that $F_T$ is a {\em square} system.

Now, suppose that $S\subset(\CC^{n_1}\times\cdots\times\CC^{n_d})^r$ consists of the known decompositions of $T$.
For a loop $\tau:[0,1]\rightarrow\CC^{n_1\cdots n_d}$
with $\tau(0) = \tau(1) = T$, consider the homotopy
$$H(v_1^1,\dots,v_d^r,s) = F_{\tau(s)}(v_1^1,\dots,v_d^r) = 0.$$
The loop $\tau$ is selected so that the solution
paths starting at the points in $S$ when $s = 0$
are nonsingular for $s\in[0,1]$.  This is the generic
behavior for paths $\tau$ since the singular locus
is a complex codimension one condition while
we are tracking along a real one-dimensional arc
$\tau(s)$ for $0\leq s\leq 1$.  The endpoints,
namely at $s =0$ and $s = 1$, of these solution paths
form a decomposition of $T$.  If a new decomposition
is found, it is added to~$S$.  The process is repeated for
many loops $\tau$.  We leave many details about path tracking
to \cite{Bertini,SWbook}.

Since $F_T$ and the homotopy $H$ is naturally
invariant under the action of the symmetric group on $r$ elements,
we only need to track one path starting from one point
from each orbit.  Each loop is usually computationally
inexpensive so we can repeat this computation many times.
Experience has shown that randomly selected loops are
typically successful at generating the requisite monodromy action
needed to obtain all decompositions starting from a single one
in a relatively small number of loops.

In the subsequent subsections, when an exact value is
reported, this means that at least~$50$ additional randomly selected
loops failed to yield any new decompositions.  Thus,
we expect that these values are sharp.
When a lower bound is reported, this means that
we have terminated the computation with
the last loop generating many new decompositions.
Thus, these lower bounds are probably quite far from
being sharp, but do show nonuniqueness.

\subsection{{The number of decompositions for perfect format tensors}}
\begin{computation}
{Tables~\ref{table3.2a} and \ref{table3.2b} summarize the results of our numerical computations which determine (with high confidence) the generic ranks and numbers of decompositions for general tensors  $d\ge 3$ satisfying \eqref{balanc:eq} with \mbox{$\prod_{i=1}^{d}n_i\le 100$}. Table~\ref{table3.3} records our results for symmetric tensors.
}
\end{computation}
\begin{table}
\[
\begin{array}{|r|r|r|}
\hline
(n_1, n_2, n_3)&\textrm{gen. rank}&\textrm{\# of decomp. of general tensor}\\
\hline
(3, 4, 5)&6&{\bf 1}\,\,\\
\hline
(3, 6, 7)&9&38\,\,\\
\hline
(4, 4, 6)&8&62\,\,\\
\hline
(4, 5, 7)&10&\geq \hbox{222,556}\,\,\\
\hline
\end{array}\]
\caption{{Results of our numerical computation for all perfect format $3$-tensors satisfying~\eqref{balanc:eq} with $\prod_{i=1}^{3}n_i\le 150$ and the number of decompositions for the generic tensor}.}\label{table3.2a}
\end{table}
\noindent

\begin{table}
\[
\begin{array}{|r|r|r|}
\hline
(n_1,\ldots,  n_d)&\textrm{gen. rank}&\textrm{\# of decomp. of general tensor}\\
\hline
(2, 2, 2, 3)& 4& {\bf 1}\,\,\\
\hline
(2, 2, 3, 4)& 6&4\,\,\\
\hline
(2, 2, 4, 5)& 8&68\,\,\\
\hline
(2, 3, 3, 4)& 8&471\,\,\\
\hline
(2, 3, 3, 5)& 9&7225\,\,\\
\hline
(3, 3, 3, 3)& 9&\hbox{20,596}\,\,\\
\hline
(2, 2, 2, 2, 4)& 8&447\,\,\\
\hline
(2, 2, 2, 3, 3)& 9& \hbox{18,854}\,\, \\
\hline
(2, 2, 2, 2, 2, 3)& 12& \geq \hbox{238,879}\,\,\\
\hline\end{array}\]
\caption{{Results of our numerical computation for all perfect format tensors with $d\ge 4$ satisfying \eqref{balanc:eq} with \mbox{$\prod_{i=1}^{d}n_i\le 100$}.and the number of decompositions for the generic tensor.}}\label{table3.2b}
\end{table}

{The generic rank
is known to be equal to the expected one}
for formats $(n,n,n)$~\cite{Lick}, which is not perfect for $n\ge 3$,
and $(2,\ldots,2)$ for at least $k\ge 5$ factors~\cite{CGG},
which is perfect if $k+1$ is a power of $2$.
A numerical check for $k=7$ shows it is~not~identifiable.

\subsection{{The number of decompositions for symmetric tensors}}\label{subsec:symmetric}

We highlight a few cases for computing
the number of decompositions of symmetric tensors.

\begin{table}
\renewcommand{\arraystretch}{1.17}
\[
\begin{array}{|r|r|r|}
\hline
\textrm{Tensor space}&\textrm{gen. rank}&\textrm{\# of decomp. of general tensor}\\
\hline
\Sym^{10}\CC^3 & 22 & 320\,\,\\
\hline
\Sym^{11}\CC^3 & 26 & 2016\,\, \\
\hline
\Sym^{5}\CC^4 & 14 & 101\,\,\\
\hline
\Sym^{3}\CC^7 & 12 & 98\,\,\\
\hline\end{array}\]
\renewcommand{\arraystretch}{1}
\caption{{Results of our numerical computation for the number of decompositions of generic tensors for some symmetric tensor formats.}}\label{table3.3}
\end{table}

In the cases of this table, Theorem 1.1 in \cite{Mella09} and the recent \cite{GM}
imply that the number of decompositions of
a general tensor is at least $2$.

For $\Sym^{d}\CC^3$ and $d=1$ or $d=2$ modulo $3$, the expectation stated in
Remark \ref{rem:tangproj} is that the number of decompositions
is divisible by $\frac{(d-2)(d-4)}{3}$.
This is confirmed for $d = 10$ with $320 = 20\cdot 16$
and $d = 11$ with $2016 = 96\cdot 21$.

\section{Pseudowitness sets and verification}\label{sec:Projection}

The approach discussed in Section~\ref{sec:Homotopy} uses
random monodromy loops to attempt to generate new decompositions.
Clearly, when showing that a format is not identifiable, one simply
needs to generate some other decomposition.
We can use the numerical approximations to generate a proof
that it is not identifiable in the perfect case using, for example,
{\tt alphaCertified}\cite{alphaCertified,alphaCertifiedSoftware}.
However, to determine the precise number of decompositions,
we simply run many monodromy loops and observe when
the number of decompositions computed stabilize.
In this section, we describe one approach for validating
the number of decompositions and demonstrate this approach
in Section~\ref{sec:366} for counting
the number of decompositions for a general
tensor of format $(3,6,6)$ of rank $8$.

\subsection{Using pseudowitness sets}\label{sec:Pseudo}

For demonstration purposes, consider counting the number
of decompositions of a general tensor of
format $(n_1,\dots,n_d)$ of rank $r$.
Consider the following where we have
removed the trivial degrees of freedom by selecting elements~to~be~$1$:
$$G := \left\{(T,v_1^1,\dots,v_d^1,\dots,v_1^r,\dots,v_d^r)~\left|~
\begin{array}{ll}
T=\sum_{i=1}^r v_1^i\otimes\cdots\otimes v_d^i, \\
(v_j^i)_1 = 1 \hbox{~for~} i = 1,\dots,r \hbox{~and~}j = 1,\dots,d-1 \end{array}
\right\}\right..$$
The graph $G$ is clearly an irreducible variety.
Hence, the image $\overline{\pi(G)}$ is also irreducible
where $\pi(T,v_1^1,\dots,v_d^1,\dots,v_1^r,\dots,v_d^r) = T$.
If $\dim G = \dim \overline{\pi(G)}$, then we know that a general
tensor~$T$ of format $(n_1,\dots,n_d)$ of rank $r$
has finitely many decompositions, namely
$$\frac{\left|\pi^{-1}(T)\cap G\right|}{r!}.$$
In particular, $\left|\pi^{-1}(T)\cap G\right|$ is the degree of a
general fiber of $\pi$ with respect to $G$ and
the denominator $r!$ accounts for the {natural} action of the symmetric group
on $r$ elements.

Using numerical algebraic geometry, computations on $G$
will be performed using a {\em witness set}, e.g., see \cite{Bertini,SWbook},
and on $\overline{\pi(G)}$ using a {\em pseudowitness set}
\cite{HW13,HW10}.  A byproduct of computing a pseudowitness set
for $\overline{\pi(G)}$ is the degree of a general fiber.

Suppose that $X\subset\CC^N$ is an irreducible variety
of dimension $k$ and $f$ is a system of polynomials in $N$ variables
such that $X$ is an irreducible component of the solution
set defined by $f = 0$.  Then, a witness set for $X$ is
the triple $\{f,\sL,V\}$ where
$\sL\subset\CC^N$ is a general linear space of codimension $k$
and $V = X\cap\sL$.  Here, ``general linear space'' means
that $\sL$ intersects the smooth points of the reduction of 
$X$ transversely so that $|V| = \deg X$.

To focus only on the case of interest, we assume
that $\pi:\CC^N\rightarrow\CC^m$ is the projection map
onto the first $m$ coordinates such that $\pi$ is
generically $\ell$-to-one on $X$,
i.e., $\dim X = \dim\overline{\pi(X)}$.
Then, a pseudowitness set for $\overline{\pi(X)}$ is the quadruple
$\{f,\pi,\sM,W\}$ where $W = X\cap\sM$ and $\sM  = \sM_{\pi}\times\CC^{N-m}$
with $\sM_{\pi}\subset\CC^m$ being a general linear space of codimension $k$.
Here, ``general linear space''
means that $\sM_{\pi}$ intersects $\overline{\pi(X)}$ transversely,
so that $|\pi(W)| = |\overline{\pi(X)}\cap\sM_{\pi}| = \deg \overline{\pi(X)}$,
and, for each $T\in \overline{\pi(X)}\cap\sM_{\pi}$, $|\pi^{-1}(T)\cap X| = \ell$.
Therefore, $\pi$ is an $\ell$-to-one map on $W$ with $|W| = \ell\cdot \deg\overline{\pi(X)}$.

Returning to the problem at hand, since we aim to compute a
pseudowitness set for~$\overline{\pi(G)}$,
we can simplify this computation by only considering the fiber
over a general curve section of~$\overline{\pi(G)}$.
If $k = \dim\overline{\pi(G)}$,
let $\sA\subset\CC^{n_1}\otimes\cdots\otimes\CC^{n_d}$
be a general linear space of codimension~\mbox{$k-1$}.
Then, $C = \pi^{-1}(\overline{\pi(G)}\cap\sA)\cap G$
and $\overline{\pi(C)} = \overline{\pi(G)}\cap\sA$
are both irreducible curves, e.g., see \cite[Thm.~3.42]{SS85}
and \cite[Thm.~13.2.1]{SWbook}, respectively.
Moreover, $\deg \overline{\pi(G)} = \deg \overline{\pi(C)}$
and the degree of a general fiber for $C$ and $G$ with respect
to~$\pi$~are~equal.

With this, we now aim to compute a pseudowitness set for
the curve~$\overline{\pi(C)}$.  
{Since $C$ is a curve, 
we compute a pseudowitness set for $\overline{\pi(C)}$ 
by intersecting $C$ with a hyperplane of the form $\pi^{-1}(\sH)$
where $\sH\subset\CC^{n_1}\otimes\cdots\otimes\CC^{n_d}$
is a general hyperplane.  Since $\pi^{-1}(\sH)$ is
invariant under the natural action of the symmetric group on $r$ elements,
we will first consider the intersection of $C$ with a general element 
of an irreducible family of hypersurfaces
that are invariant under the same action and contains 
hyperplanes of the form $\pi^{-1}(\sH)$.
This is sufficient by the results of \cite{CoeffParam} and 
simplifies the computation since only one point in each
orbit needs to be computed, i.e., only one point 
out of every $r!$ based on the natural action of the symmetry
group on $r$ elements.
} 

\begin{remark} 
{
At first sight, by reducing down to the curve case,
it might seem that some technicalities could be avoided by 
choosing a general hyperplane among the hyperplanes invariant under 
the symmetric group and only needing to compute one point in each orbit.   
}
However, there could be some difficulties 
that arise with this.

The first difficulty is that the set of hyperplanes invariant under a finite group does not have to be irreducible.  
{
For example,} let $G\subset \CC^*$ denote the sixth roots of unity.
Consider the action of $G$ on $\CC^2$ {
defined by $g\cdot(z,w)\mapsto (g^2z,g^3w)$. }
In this case, the set of invariant hyperplanes consists of 
two points, i.e., the $z$ axis and the $w$ axis. 

{
Even when there is a large family of invariant hyperplanes that make sense, 
the second difficulty is that none of them need be invariant enough to 
intersect the algebraic set of interest transversely in the degree number of points.
}
For example, let $X\subset \CC^2$ be the solution set of $z^2-w^3$
{which is a curve of degree $3$}.
Let $G$ denote $\{1,-1\}$ 
{
which acts on $\CC^2$ by $g\cdot (z,w)\mapsto (gz,w)$.}
Note that $X$ is invariant under $G$.  The $G$-invariant hyperplanes consist of two components.
One component is made up of the fibers of the 
{projection} map $(z,w)\mapsto w$
{which meet} $X$ in two points.
The other component consists of the $w$ axis which is not transversal to $X$. 

To avoid these potential difficulties, we first define an irreducible family
of hyperplanes that contains the invariant hyperplanes of interest.
{Then, the} 
general theory in \cite{CoeffParam} shows that this process 
{computes} the requisite points needed in the 
construction of a pseudowitness~set.
\end{remark}

\subsection{Tensors of format $(3, 6, 6)$ of rank $8$.}\label{sec:366}

The tensors of format $(3, 6, 6)$ have generic rank $9$
in which a general tensor of this format has
infinitely many decompositions.
In \cite{CMO}, the open problem of computing the number of
tensor decompositions of a general tensor of rank $8$
of format $(3,6,6)$ was formulated.
To the best of our knowledge, this is probably the last open case
when a generic tensor of some rank strictly
smaller than the generic one is not identifiable.
Theorem 3.5 of \cite{CMO} proved that the
number of decompositions is $\ge 6$.
Moreover, \cite{CMO} showed that the number of
decompositions of format $(3,6,6)$ of rank $8$ is
equal to the number of decompositions of a general tensor
in $\Sym^3\CC^3\otimes\CC^2\otimes\CC^2$, which is perfect
with generic rank~$8$.  We use the approach from
\S\ref{sec:Pseudo} to show exactly $6$ decompositions.
Another approach based on a multihomogeneous trace test
in \cite{HR} confirms this result.

To use the approach
from \S\ref{sec:Pseudo}, consider the Veronese embedding
$$v_3:\CC^3\rightarrow\CC^{10} \hbox{~~where~~}
(x,y,z) \mapsto (x^3,x^2y,x^2z,xy^2,xyz,xz^2,y^3,y^2z,yz^2,z^3).$$
We picked a random line $\sL\subset\CC^{40}$ and
consider the irreducible curve
{\footnotesize
$$
C := \left\{(T,a_1,b_1,c_1,\dots,a_8,b_8,c_8)~\left|~
\begin{array}{l}T = \sum_{i=1}^8 v_3(a_i)\otimes b_i\otimes c_i \in \sL
\\ (a_i)_1 = (b_i)_1 = 1 \hbox{~for~} i = 1,\dots,8 \end{array}
\right\}\right.
\subset \sL\times \left(\CC^3\times\CC^2\times\CC^2\right)^8.$$
}

{To compute a pseudowitness set for $\overline{\pi(C)}$, we
need to compute $C\cap(\sM_\pi\times(\CC^3\times\CC^2\times\CC^2)^8)$
where $\sM_\pi\subset\CC^{40}$ is a general hyperplane.
Consider the irreducible family of hyperplanes 
defined by the vanishing of linear equations of the form
\begin{equation}\label{eq:Hyperplanes}
\sum_{j=1}^{40} \alpha_j T_j + \sum_{j=1}^3 \beta_j\sum_{k=1}^8 (a_k)_j
+ \sum_{j=1}^2 \gamma_j\sum_{k=1}^8 (b_k)_j + \sum_{j=1}^2 \delta_j\sum_{k=1}^8 (c_k)_j = \epsilon.
\end{equation}
By construction, this family is invariant under the natural action of 
the symmetric group $S_8$ and contains all hyperplanes
of the form $\sM_\pi\times(\CC^3\times\CC^2\times\CC^2)^8$.
After picking a random hyperplane $\sH$
of the form \eqref{eq:Hyperplanes}
and starting with one point on $C\cap\sH$,}
we used monodromy loops via {\tt Bertini} to compute additional points
in $C\cap\sH$ which stabilized to $1020\cdot 8! = \hbox{41,126,400}$.
The trace test~\cite{trace} confirms that this set
of points is indeed equal to $C\cap\sH$.
{After selecting a random hyperplane $\sM_\pi$, 
we then computed 
$W = C\cap(\sM_\pi\times(\CC^3\times\CC^2\times\CC^2)^8)$
by deforming from $C\cap\sH$ which yielded
$|W| = 6\cdot 8! = \hbox{241,920}$.
Since $\overline{\pi(C)}\cap\sM_\pi = \sL$ which has degree~$1$,
$|\pi(W)| = 1$ thereby showing} exactly $6$ decompositions.

{
We summarize the result of this computation.

\begin{computation}
Numerical algebraic geometry together with randomly selected
linear spaces show {\it (with high confidence)} that a general tensor 
in $\Sym^3\CC^3\otimes\CC^2\otimes\CC^2$ has exactly $6$ decompositions.
\end{computation}
}

\section{Tensor decomposition via apolarity}\label{sec:Eigen}

In \cite{OeOtt}, a technique generalizing Sylvester's algorithm was implemented
by considering the kernel of the catalecticant map, which in turn is a graded
summand of the apolar ideal.
In principle, this apolarity technique can be used for any algebraic variety.

\subsection{A uniform view of (Koszul) flattenings}
Let $V$ {and} $V_{i}$ be arbitrary finite dimensional vector spaces over $\CC$ of dimensions $n$ {and} $n_{i}$, respectively.
For $-n\leq p\leq n$, let $\bw{p}V$ denote fundamental representations of $GL(V)$
where we interpret $\bw{p} = \bw{-p} V^{*}$ when $p<0$.
For a multi-index $I\in \ZZ^{d}$, let $V_{I}$ denote the tensor product of fundamental representations
\[
V_{I}:= \bw{i_{1}} V_{1} \o \bw{i_{2}} V_{2} \o \dots \o \bw{i_{d}} V_{d}
.\]
Note that $V_{1^{d}}:=V_{1}\o \dots \o V_{d}$.
We may assume, up to reordering, that $i_j\ge 0$ for $j=1,\ldots, h$,
$i_j<0$ for $j=h+1,\ldots, d$.
We obtain linear maps $K_{p}\colon \bw{p}V \to \bw{p+1} V$
that depend linearly on $V$
by way of the Koszul complex.  Specifically, for $v \in V$ and $\phi \in \bw{p}V$ define
\[
K_{p}(v)(\phi) = \phi \wedge v\textrm{\ for\ }p\ge 0,
\]
\[
K_{p}(v)(\phi) = \phi(v)\textrm{\ for\ }p<0.
\]

Now we consider the tensor product of many Koszul maps, which are linear maps on tensor products of fundamental representations  that depend linearly on  $T\in V_{(1,\dots,1)}:=V_{1}\o\dots \o V_{d}$:
\[
K_{I}(T)\colon V_{I} \to V_{I + 1^{d}}
.\]
For indecomposable elements $v_{1}\o\dots\o v_{d} \in V_{1}\o\dots \o V_{d}$ and $\phi_{1}\o \dots \o \phi_{d} \in V_{I}$~define
\begin{equation}
K_{I}(v_{1}\o\dots\o v_{d})(\phi_{1}\o \dots \o \phi_{d}) = \bigotimes_{j=1}^h(\phi_{j}\wedge v_{j}) \o \bigotimes_{j=h+1}^d (\phi_{j}( v_{j}))
.\label{eq:defKI}\end{equation}
The definition of $K_{I}$ is extended by bi-linearity. {We often drop the argument $(v)$ or $(T)$ of $K_{p}$ or $K_{I}$ when the tensor to be flattened is not specified and the linear dependence is understood.}
From this definition it is clear that the image of $K_{I}(v_{1}\o\dots\o v_{d})$ is isomorphic~to
\[
\bigotimes_{j=1}^h \left(\bw{i_{j}} \left(V_{j}/\langle v_{j} \rangle \right) \o (v_{j})\right)\o \bigotimes_{j=h+1}^d \left(\bw{-i_{j}-1} \left(v_{j}^\perp \right)\right).
\]
A consequence of dimension counting, bi-linearity of $K_{I}$, and
sub-additivity of matrix rank is the following, which
is essentially described in \cite[Prop.~4.1]{LO11}.
\begin{prop0}\label{prop:LanOtt}
Suppose $T \in V_{1,\dots,1}$ has tensor rank $r$. Let $i_j\ge 0$ for $j=1,\ldots, h$,
$i_j<0$ for $j=h+1,\ldots, d$. Then the Koszul flattening $K_{I}(T)\colon V_{I} \to V_{I+1^{d}}$ has rank at most
\[
r_{I}:=r\cdot \prod_{j=1}^{h} \binom{n_{j}-1}{i_{j}} \cdot \prod_{j=h+1}^{d} \binom{n_{j}-1}{-i_{j}-1}
.\]
In particular,  the $(r_{I}+1) \times (r_{I}+1)$ minors of $K_{I}(T)$ vanish.
This is meaningful provided that $r_{I}<\min\{\dim V_{I}, \dim  V_{I+1^{d}}\}$.
\end{prop0}
{Let $m = \prod_{j=1}^{h} \binom{n_{j}-1}{i_{j}} \cdot \prod_{j=h+1}^{d} \binom{n_{j}-1}{-i_{j}-1}$. The proposition says that for all tensors $T$ with $\rank(T)=r$ we have $\rank(K_{I}(T) ) \geq r \cdot m$.
Thus, Koszul flattenings potentially provide the most useful information whenever the following ratio is maximized:
\[
\min\left\{
\dim V_{I}, \dim  V_{I+1^{d}}
\right\}
\;/\; m
.\]
When $rm \leq \min\left\{
\dim V_{I}, \dim  V_{I+1^{d}} \right\}$ we say that the flattening can detect rank $r$ since we expect that the flattening $K_{I}(T)$ will have different ranks when $T$ has rank either $r-1$ or $r$ {(this is true in all the examples we have tried, including the ones occurring in Theorems \ref{345:thm} and
\ref{thm2223})} and there is no obstruction to this based solely on the size of the flattening matrix.
}

\subsection{Apolarity Lemma for Koszul flattenings}
Recall from (\ref{eq:defKI})
that{
\[
 K_{I}(v_{1}\o \dots \o v_{d}) \equiv 0 \Longleftrightarrow \bigotimes_{j=1}^h (\phi_{j}\wedge v_{j})\o  \bigotimes_{j=h+1}^d (\phi_{j}(v_{j}))= 0\] }
for all {pure tensors} $\phi \in V_{I}$.

It is useful to look at tensors in the kernel of $K_{I}(T)$ as linear maps. With this aim, we need to distinguish the negative and nonnegative parts of $I \in \ZZ^{d}$.  So let $N \sqcup P =\{1,\dots,d\}$ be the set partition such that $-I_{N}\in \ZZ_{> 0}^{d}$,
$I_{P} \in \ZZ_{\ge 0}^{d}$ and the notation $I_{P}$ (resp. $I_N$) is the vector in $\ZZ^{d}$ gotten by keeping the elements of $I$ in the positions $P$ (resp. $N$) and {setting the rest of the entries to zero}. We also let $1_{P}$ denote the vector with ones in the positions denoted by the index $P$ (and zero elsewhere), and similarly for $1_{N}$.
 With this, we may identify
$
V_{I} = V_{I_{P} + I_{N}}= \Hom (V_{-I_{N}},V_{I_{P}}),
$
and consider the Koszul flattening of $T\in V_{(1,\dots, 1)}$ as
\[K_{I}(T) \colon \Hom (V_{-I_{N}},V_{I_{P}}) \to  \Hom(V_{-I_{N}+1_{N}},  V_{I_{P}+1_{P}}) .\]
{Up to reordering the factors, $K_{I}(T)$ is defined on decomposable elements $\left(\bigotimes_{j\in I_N}w_j\right)$ by }

\[K_{I}\left(v_{1}\o\dots\o v_{d}\right)(\psi)\left(\bigotimes_{j\in I_N}w_j\right)=\psi\left(\bigotimes_{j\in I_N}(w_j\wedge v_j)\right)\bigwedge\left( \bigotimes_{j\in I_p}v_{j}\right) \quad
\forall \psi\in \Hom (V_{-I_{N}},V_{I_{P}}).\]

In our setting, \cite[Prop.~5.4.1]{LO11} yields the following lemma
(see (\ref{apol:eq}) for a concrete case).
Since this technique refers to a vector bundle, it could be called ``nonabelian'' apolarity,
in contrast with classical apolarity which refers to a line bundle (see \cite[Ex.~5.1.2]{LO11} and \cite[\S~4]{OeOtt}).
\begin{lemma0}[Apolarity Lemma]\label{apolarity:lem} Suppose $T = \sum_{s=1}^{r}v_{1}^{s}\o \dots \o v_{d}^{s}$.
\[
\mbox{\small $
\ker K_{I}(T) \supset  \{\psi\in \Hom (V_{-I_{N}},V_{I_{P}}) \mid
 \psi\left( V_{-I_N+1_N}\wedge\bigotimes_{j\in N} v_{j}^{s}\right) \wedge \left(\bigotimes_{j \in P} v_{j}^{s}\right)
=0\textrm{\ for\ }s=1,\ldots, r\}.$}
\]

\end{lemma0}
\begin{proof}
{Pick $\psi$ in the space on the right hand side of the inclusion.} Choosing any $w_j\in V_{-I_N+1_N}$ for every $j\in I_N$, we have
$$K_{I}\left(\sum_{s=1}^{r}v_{1}^{s}\o \dots \o v_{d}^{s}\right)(\psi)\left(\bigotimes_{j\in I_N}w_j\right) =
\sum_{s=1}^{r}
\psi\left(\bigotimes_{j\in I_N}(w_j\wedge v_j^s)\right)\bigwedge\left( \bigotimes_{j\in I_p}v_{j}^s\right)$$
and each summand vanishes by the assumption.\end{proof}

\subsection{The $3\times 4\times 5$ case}\label{Sec:345}

Let us denote the three factors as
$A=\CC^3$, $B=\CC^4$, $C=\CC^5$.
The following are all possible {relevant} non-redundant Koszul flattenings (up to transpose), {which all depend linearly on $A\o B\o C$.}
\[
K_{(0,-1,-1)}\colon (B\o C)^{*}  \to  A , \quad
K_{(-1,0,-1)} \colon (A\o C)^{*}  \to B ,\quad
K_{(-1,-1,0)}\colon (A\o B)^{*}  \to C
\]
\[
{K_{(1,0,-1)}
\colon C^{*}\o A \to B\o \bw{2} A  \quad \text{,}\quad
K_{(-1,0,2)}
\colon A^{*}\o \bw{2}C \to B\o \bw{3} C}
\]
\[
K_{(0,1,-1)}
\colon C^{*}\o B \to A\o \bw{2} B \quad \text{,}\quad
K_{(-1,1,0)}
\colon A^{*}\o B \to C\o \bw{2} B,
\]
\[
K_{(-1,0,1)}
\colon A^{*}\o C \to B\o \bw{2} C \quad \text{,}\quad
K_{(0,-1,1)}
\colon B^{*}\o C \to A\o \bw{2} C,
\]

Consider
\[
K_{u} (T) \colon \bw{u_{1}}A\o \bw{u_{2}}B \o \bw{u_{3}}C \to  \bw{u_{1}+1}A\o \bw{u_{2}+1}B \o \bw{u_{3}+1}C
,\]
where, for any vector space $V$, we interpret negative exterior powers by asserting
that \mbox{$\bw{s} V = \bw{-s} V^{*}$} if $s<0$.

For example $K_{0,-1,-1}(a\o b\o c)$ has image
\[   (\bw{0}A \wedge a)\o  (B^*( b)) \o (C^*( c)) \subset \bw{1}A\o \bw{0}B \o \bw{0}C.
\]
The factor  $ (B^*( b)) \o (C^*( c))$ is {just any scalar, obtained
by contracting $b$ with $B^*$, and $c$ with~$C^*$.} We are left with $(\bw{0}A \wedge a) = \langle a \rangle$, which is 1-dimensional.

Another example is $K_{0,1,-1}(a\o b\o c)$ which has image
\[   (\bw{0}A \wedge a)\o  (\bw{1}B\wedge b) \o (C^*( c)) \subset \bw{1}A\o \bw{2}B \o \bw{0}C.
\]
The factor  $C^*( c)$ is a scalar obtained by contracting $c$ with $C^*$.
{This leaves} $(\bw{0}A \wedge a) = \langle a \rangle$ tensored with $(\bw{1}B\wedge b) \subset  \bw{2}B$, but $(\bw{1}B\wedge b) \cong (B/b) \o \langle b \rangle $, which is 3-dimensional.

In general, the image of $K_{u}(a\o b \o c)$ has dimension
\begin{equation}\label{eq:image}
\binom{\dim A -1}{f(u_{1})} \binom{\dim B -1}{f(u_{2})} \binom{\dim C -1}{f(u_{3})}
.\end{equation}
where $f(x)=\left\{\begin{array}{cc}x&\textrm{\ if\ }x\ge 0\\ -x-1&\textrm{\ if\ }x< 0\end{array}\right.$.
On the other hand, the maximum rank that $K_{u}$ can have is the minimum of the dimensions of the source and the target, or
\begin{equation}\label{eq:max}
\min\left\{\binom{\dim A}{|u_{1}|} \binom{\dim B }{|u_{2}|} \binom{\dim C }{|u_{3}|},
\binom{\dim A }{|u_{1}+1|} \binom{\dim B }{|u_{2}+1|} \binom{\dim C }{|u_{3}+1|} \right\}.
\end{equation}
{Since the matrix rank of a Koszul flattening of a tensor $T$ is bounded above by a multiplicative factor \eqref{eq:image} of the tensor rank of $T$,
 the maximum tensor rank that a Koszul flattening can detect by a drop in matrix rank} is the ratio of \eqref{eq:max} and~\eqref{eq:image}.
For convenience we record the dimensions and the multiplication factor (\ref{eq:image}) for each flattening.

\begin{center}
\begin{tabular}{|c|c|c|c|}\hline
map & size & mult-factor & max tensor rank detected\\
\hline
$K_{(0,-1,-1)} $& $ 3\times 20 $ & 1 & 3\\\hline
$K_{(-1,0,-1)} $& $ 4\times 15 $ & 1 & 4\\\hline
$K_{(-1,-1,0)} $& $ 5\times 12 $ & 1 & 5\\\hline
$K_{(1,-1,0)} $& $ 15\times 12 $ & 2 & 6\\\hline
$K_{(1,0,-1)} $& $ 12\times 15 $ & 2 & 6\\\hline
$K_{(0,1,-1)} $& $ 18\times 20 $ & 3 & 6\\ \hline
$K_{(-1,1,0)} $& $ 12\times 30 $  & 3 & 4\\\hline
$K_{(-1,0,1)} $& $ 40\times 15 $& 4 & 4\\\hline
$K_{(0,-1,1)} $& $ 30\times 20 $ & 4 & 5\\\hline
$K_{(-1,0,2)} $& $ 40\times 30 $ & 6 & 5\\\hline
$K_{(0,-1,2)} $& $ 30\times 40 $ & 6 & 5\\\hline
$K_{(0,-1,2)} $& $ 30\times 40 $ & 6 & 5\\\hline
\end{tabular}\end{center}\bigskip

We see that the only maps that {might} distinguish between tensor rank 5 and 6 are $K_{(1,-1,0)}$, $K_{(1,0,-1)}$, and $K_{(0,1,-1)}$.
Since $\bw{2}A \cong A^{*}$, the first two maps are transposes of each other:
\[K_{(1,-1,0)}  = (K_{(1,0,-1)})^{t}.\]
Thus, we proceed by considering $K_{(1,0,-1)}$ and $K_{(0,1,-1)}$.

In our case, Apolarity Lemma \ref{apolarity:lem} says that
\begin{equation}
\label{apol:eq}\ker K_{1,0,-1}( \textstyle \sum_{i=1}^{s} a_ib_ic_i)\supset\{\phi\in Hom(C,A)|\phi(c_i)\wedge a_i=0\textrm{\ for\ }i=1,\ldots, s\}.\end{equation}
and
\[\ker K_{0,1,-1}( \textstyle \sum_{i=1}^{s} a_ib_ic_i)\supset\{\phi\in Hom(C,B)|\phi(c_i)\wedge b_i=0\textrm{\ for\ }i=1,\ldots, s\}.\]
Equality should hold for honest decompositions, see \cite[Prop.~5.4.1]{LO11}.

With this setup, we now present a proof of Theorem~\ref{345:thm}.

\begin{proof}[Proof of Theorem \ref{345:thm}]
For general $f\in A\otimes B\otimes C$, $K_{1,0,-1}(f)$ is surjective and $\ker K_{1,0,-1}(f)$ has dimension  $\dim Hom(C,A)-\dim \wedge^{2}A\otimes B=
15-12=3$. To complete the proof we  interpret the linear map $K_{1,0,-1}(f)$ as a map between sections of vector bundles.
Let $X=\PP(A)\times\PP(B)\times\PP(C)$, endowed with the three projections $\pi_A$, $\pi_B$, $\pi_C$ on the three factors.
We denote  $\O(\alpha,\beta,\gamma)=\pi_A^*\O(\alpha)\otimes\pi_B^*\O(\beta)\otimes\pi_C^*\O(\gamma) $.
Let $Q_A$ be the pullback of the quotient bundle on $\PP(A)$.

Let $ E=Q_A\otimes\O(0,0,1)$ and $L=\O(1,1,1)$. Note that $E$ is a rank two bundle
on~$X$.
As in~\cite{LO11,OeOtt}, the map $K_{1,0,-1}(f)$ can be identified with the contraction
$$K_{1,0,-1}(f)\colon H^0(E )\rig{} H^0(E^*\otimes L)^{*}$$
which depends linearly on {$f\in H^{0}(L)^*$}.

The general element in $H^0(E)$ vanishes on a codimension two subvariety of $X$
which has the homology class $c_2(E)\in H^*(X,\ZZ)$.
The ring $H^*(X,\ZZ)$ has three canonical generators $t_A, t_B, t_C$ and it can be identified with
$\ZZ[t_A,t_B,t_C]/(t_A^3, t_B^4, t_C^5)$. {Since the rank of $Q_A$ is $2$
and $\O(0,0,1)$ is a line bundle, we have the well known formula
 $c_2(E)=c_2(Q_A)+c_1(Q_A)c_1(\O(0,0,1))+c_1^2(\O(0,0,1))$.
We have $c_1(Q_A)=t_A$, $c_2(Q_A)=t_A^2$, $c_1(\O(0,0,1))=t_C$ and we compute }
$c_2(E)=t_A^2+t_At_C+t_C^2$. Hence three general sections of $H^0(E)$
have their common base locus given by $c_2(E)^3=\left(t_A^2+t_At_C+t_C^2\right)^3=6t_A^2t_C^4$.
This coefficient $6$ coincides with the generic rank and it is the key of the computation. {Terracini's lemma {(see \cite[Cor.~5.3.1.2]{Land12})}, using
a simple tangent space computation,
verifies that the generic rank is indeed 6.
We pick a tensor $f=\sum_{i=1}^6a_ib_ic_i$ constructed with six random points
$(a_i, b_i, c_i)\in\PP(A)\times\PP(B)\times\PP(C)$ for $i=1,\ldots, 6$.
By using {\tt Macaulay2} (see the {\tt Macaulay2} file attached at the arXiv {version of this paper}) we may compute that the variety $$\{(a,c)\in \PP(A)\times\PP(C)|
\phi(c)\wedge a=0\quad\forall \phi\in\ker K_{1,0,-1}(f)\}$$ consists of
the union of six points $(a_i, c_i)$ for $i=1,\ldots, 6$, and this holds even scheme-theoretically. It follows that equality holds in (\ref{apol:eq}) with $s=6$, { indeed $K_{1,0,-1}(f)$ is surjective and has $3$-dimensional kernel, while the right-hand side of (\ref{apol:eq}) has dimension 
at least $\dim Hom(C,A)-6\cdot 2 = 3$ .} 
Moreover the common base locus of
$\ker K_{1,0,-1}(f)$ is given by~$6$~linear spaces  $\{a_i\}\times\PP(B)
\times\{ c_i\}$ for $i=1,\ldots, 6$}.
 
 { We claim that} the common base locus of
$\ker K_{1,0,-1}(f)$ is given by $6$ { linear spaces as above} for general tensor $f$. { Indeed, the common zero locus of three sections of $E$ can be seen as the zero locus of a section of $E^{\oplus 3}$. The dimension of the zero locus of a section of $E^{\oplus 3}$  is at least $9-6=3$ and it is upper semicontinous with respect to the section, hence the common base locus of
$\ker K_{1,0,-1}(f)$ is a pure-dimensional $3$-fold for general $f$.
Since the top Chern class of $E^{\oplus 3}$ is $c_2(E)^3$,
which we computed to be $6t_A^2t_C^4$, we know that the common base locus of
$\ker K_{1,0,-1}(f)$ comes from a $0$-dimensional scheme of degree $6$ on $\PP(A)\times\PP(C)$,
after a pullback with $\pi_B$, by \cite[Prop. 14.1 (b)]{Ful}. Since this $0$-dimensional scheme consists of $6$ distinct points for the tensor $f$ that we have picked, this property remains true for general~$f$.
}

In particular, the decomposition {$f=\sum_{i=1}^6a_i\otimes b_i\otimes c_i$}
has a unique solution (up to scale) for $a_i$ {and} $c_i$.
{After $a_i$, $c_i$ have been determined,} the remaining vectors $b_i$ can be recovered uniquely by solving a linear system.
\end{proof}
\begin{remark}
If we attempt to repeat the same proof using $K_{0,1-1}$ in place of $K_{1,0,-1}$  most parts go through unchanged.
The map $K_{0,1-1}(T) \colon C^* \o B \to A \o \bw{2} B$
 is $18 \times 20$, and general element $T$ produces a 2-dimensional kernel.
Then, we consider the intersection of two general sections of $H^0(E) = \Hom(C,B),$ where now $E = Q_{B}\o \O(0,0,1)$.
The top Chern class of $E$ is (by a similar calculation as in the proof of Theorem \ref{345:thm})
\[4 t_B^3 t_C^3 + 3 t_B^2 t_C^4
.\]
This gives that the common base locus of $\ker K_{0,1-1}(T)$ is given by a degree 7 curve on the Segre product   $\Seg(\PP {(C)} \times \PP {(B)})$. This curve necessarily contains the 6 points needed to decompose $T$, but
information from another Koszul flattening is needed to find them.
\end{remark}

\begin{remark}
For the $(3, 4, 5)$ format, we can even decompose a general tensor $T$ of any
rank~$r$ between 1 and 6. The trick is to add {to} $T$ the sum
of $6-r$ general decomposable tensors, find the unique decomposition with the
algorithm described in the proof of Theorem~\ref{345:thm}, and
subtract the $6-r$ tensors that have to appear in the decomposition.
Unfortunately, this technique cannot
work in other cases if we do not have a tensor decomposition to start with.
\end{remark}

\subsection{The $\geq 4$ factor case}

We have seen that general tensors of format $(2, n, n)$ and
$(3, 4, 5)$ are identifiable.
We asked if there {are} other formats with this property
when there are $\geq 4$ factors.
To our surprise, the numerical homotopy method predicted an
additional case where identifiability holds.
Our construction of Koszul flattenings and
multi-factor apolarity above allows us to prove this fact, which
we prove next.

\subsection{The $2\times 2\times 2 \times 3$ case}
For this part, let $A \cong B\cong C \cong  \CC^{2}$ and $D\cong \CC^{3}$.
Because of the small dimensions we are considering, the number of interesting Koszul flattenings for tensors in $A\o B\o C\o D$ is limited to the following maps, which depend linearly on $A\o B\o C\o D$.

The 1-flattenings (and their transposes):
\[
K_{-1,0,0,0}\colon A^{*}\to B\o C \o D ,\quad
K_{0,-1,0,0}\colon B^{*}\to A\o C \o D ,\]\[
K_{0,0,-1,0}\colon C^{*}\to A\o B \o D ,\quad
K_{0,0,0,-1}\colon D^{*}\to A\o B\o C,
\]
which detect a maximum of rank 2 in the first 3 cases and a maximum of rank 3 in the last.

The 2-flattenings (and their transposes):
\[
K_{0,0,-1,-1}\colon C^{*}\o D^{*}\to A \o B,\quad
K_{0,-1,0,-1}\colon B^{*}\o D^{*}\to A \o C,\]\[
K_{-1,0,0,-1}\colon A^{*}\o D^{*}\to B \o C.
\]
The maps are all $4\times 6$ and detect a maximum of tensor rank 4.
\begin{remark}
{For (2,2,2,2), it is known} that only two of the three 2-flattenings are algebraically independent, and the dependency of the third on the other two is ``responsible'' for the defectivity of the $3^{\rm rd}$ secant variety $\sigma_{3}(\PP^{1}\times \PP^{1}\times \PP^{1}\times \PP^{1})$.
{This secant variety has one dimension} less than expected.
This type of Segre variety was indeed studied by C. Segre \cite{Segre1920}.
\end{remark}

The higher Koszul flattenings:
\[
K_{-1,0,0,1}\colon A^{*} \o D \to B \o C \o \bw{2}D,\quad
K_{0,-1,0,1}\colon B^{*} \o C \to A\o  C \o \bw{2}D,\]\[
K_{0,0,-1,1}\colon C^{*} \o D \to A\o B \o \bw{2}D
\]
These maps are all $12\times 6$, and detect a maximum of rank 3.

{We will proceed with the 2-flattenings
in the following proof of Theorem~\ref{thm2223}
since they are the only flattenings that detect the difference
between rank 3 and 4.}

\begin{proof}[Proof of Theorem~\ref{thm2223}]
Terracini's lemma {(see \cite[Cor.~5.3.1.2]{Land12})}, using
a tangent space computation, verifies the well known fact that the generic rank is $4$.
Suppose $T\in A\o B \o C \o D$ is general among tensors of rank 4 and write $T = \sum_{i=1}^{4} a_{i}\o b_{i} \o c_{i} \o d_{i}$.

First consider the case $K_{0,0,-1,-1}\colon C^{*}\o D^{*}\to A \o B$. If $T$ is general of rank $4$, then $K_{0,0,-1,-1}(T)$ has rank 4, and must have a 2-dimensional kernel. {Now, we apply \eqref{apol:eq}}. 
The points $\{c_{i}\o d_{i}\}$ must be contained in the common base locus of the elements in the kernel of $K_{0,0,-1,-1}(T)$.
Let $X=\PP(A)\times\PP(B)\times\PP(C)\times\PP(D)$, endowed with the four projections $\pi_A$, $\pi_B$, $\pi_C$, $\pi_D$ on the four factors.
We denote  $\O(\alpha,\beta,\gamma,\delta)=\pi_A^*\O(\alpha)\otimes\pi_B^*\O(\beta)\otimes\pi_C^*\O(\gamma)\otimes\pi_D^*\O(\delta) $.
Consider the line bundle $E = \Osh(0,0,1,1)$ and $L = \Osh(1,1,1,1)$  over $X$. The ring $H^*(X,\ZZ)$ has four canonical generators $t_A, t_B, t_C, t_D$ and it can be identified with
$\ZZ[t_A,t_B,t_C, t_D]/(t_A^2, t_B^2, t_C^2, t_D^3)$. 
We have $c_1(E)=t_C+t_D$. { As in the proof of Theorem~\ref{345:thm},
the common base locus of the $2$-dimensional kernel of $K_{0,0,-1,-1}(T)$
may be seen as the zero locus of a section of $E^{\oplus 2}=\Osh(0,0,1,1)^{\oplus 2}$, which has top Chern class $c_1^2(E)$.} Since $c_1^2(E)=2t_Ct_D+t_D^2$, we get
{$c_1^2(E)t_C=1 \cdot t_C t_D^2$, $c_1^2(E)t_D=2\cdot t_C t_D^2$.} 
It follows that two general sections of $E$ have common base locus given by a cubic curve, denoted $\mathcal{C}_{C,D}$ of bi-degree $(1,2)$ on $\Seg(\PP {(C)} \times \PP {(D)})$,
pulled back to a codimension $2$ subvariety of $X$.  The projection to $\PP (D)$ is a conic, which we denote $\mathcal{Q}_{C}$.

Similarly for the next 2-flattening, $K_{0,-1,0,-1}\colon B^{*}\o D^{*}\to A \o C,$ we repeat the same process, where all the dimensions and bundles are the same except for a change of roles of $C$ and $B$. By the same method we obtain another conic $\mathcal{Q}_{B}$ in {$\PP (D)$}.

Finally, if $\mathcal{Q}_{C}$ and $\mathcal{Q}_{B}$ are general,
B\'ezout's theorem implies that they intersect in 4 points in $\PP (D)$,
say {$\{[d_{1}],[d_{2}],[d_{3}],[d_{4}]\}$. Like in the proof of Theorem \ref{345:thm} we can check, using a {\tt Macaulay2} script, that starting from $T=\sum_{i=1}^4a_ib_ic_id_i$ with $(a_i, b_i, c_i, d_i)\in X$ for $i=1,\ldots, 4$, we get that the conics found by the above procedure intersect exactly in $\{[d_{1}],[d_{2}],[d_{3}],[d_{4}]\}$. } { By semicontinuity,  like in the proof of Theorem \ref{345:thm} ,
we have  that the intersection between the common base locus of the kernel of $K_{0,0,-1,-1}(T)$ and the common base locus of the kernel of $K_{0,-1,0,-1}(T)$, corresponding to the class $(2t_Ct_D+t_D^2)(2t_Bt_D+t_D^2)=4t_Bt_Ct_D^2$, consists of the pullback under $\pi_A$ of four distinct points on 
$\PP(B)\times\PP(C)\times\PP(D)$, for general $T$, namely of the $4$ linear spaces
$ \PP(A)\times \{b_{i}\} \otimes\{ c_{i}\} \times \{d_{i}\}$.} 
{ In particular, the decomposition {$T=\sum_{i=1}^4a_i\otimes b_i\otimes c_i\otimes d_i$}
has a unique solution (up to scale) for $b_i$, $c_i$ and $d_i$. After $b_i$, $c_i$, $d_i$ have been determined, the remaining vectors $a_i$ can be recovered uniquely by solving the linear system
$
T = \sum_{i=1}^{4} a_{i}\o b_{i}\o c_{i}\o d_{i}
$.}
\end{proof}

\def\null0
{
\subsection{Other formats}
In this subsection we provide evidence that the nonabelian apolarity technique with Koszul flattenings does not predict any further cases of generic identifiability of tensors in $\CC^{n_{1}\times\dots\times n_{d}}$.
For this we suppose that everything is sufficiently general.

First, we require that the tensor format be perfect, so
the fraction $\frac{\prod_{i=1}^dn_i}{1+\sum_{i=1}^d(n_i-1)}$ is an integer.

Second, we require that the Koszul flattening $K_{I}(T)$ must have a non-trivial kernel for the generic rank, and we require that the maximal rank of $K_{I}(T)$ is maximal precisely when $T$ has the generic rank (and not for sub-generic rank $T$).
\[
\min\left\{
\prod_{j=1}^{d} \binom{n_{j}}{i_{j}} ,\prod_{j=1}^{d} \binom{n_{j}}{i_{j}+1}
\right\}
\;/\; \prod_{j=1}^{d} \binom{n_{j}-1}{i_{j}}
.\]

Third, we consider the line bundle $L = (1,\dots,1)$ and the vector bundle $E_{I} = Q_{V_{1}}(-i_{1})\o \dots \o Q_{V_{d}}(-i_{d})$ associated to the Koszul map $K_{I} \colon V_{I} \to V_{I+ 1^{d}}$.
We consider the bundle map
\[
H^{0}(E) \to H^{0}(E\o L),
\]
and let $X_{I}$ denote the variety embedded by the sections of the bundle $E_{I}$.

 For $T\in V_{1^{d}}$ general, the kernel of $K_{I}(T)$ has dimension
\[M_{I}:=
\prod_{j=1}^{d} \binom{n_{j}}{i_{j}}
-\prod_{j=1}^{d} \binom{n_{j}-1}{i_{j}}
.\]
In order for the common base locus of $M_{I}$ general sections of $E_{I}$  (over $X_{I}$) to be predicted to be non-empty, we need $M_{I}*\rank E_{I} \leq \dim X_{I}$.

{Is this already enough restrictions to reduce to finitely many cases?}

After all these restrictions, we wonder if there are only finitely many cases left.  It remains to consider, like in the format $(2,2,2,3)$-case, when a combination of different flattenings must be used to compute the decomposition.
}

\section{Conclusion}\label{sec:Conclusion}

By using a numerical algebraic geometric approach based on monodromy
loops, we are able to determine the number of decompositions
of a general tensor.
Since this approach determined that general tensors of format $(3,4,5)$ and
$(2,2,2,3)$ have a unique decomposition, we developed
explicit proofs for these two special cases.
With the classically known generically identifiable case of matrix
pencils, i.e., format $(2,n,n)$, we conjecture these are the only
cases for which a general tensor has a unique decomposition.

We are currently researching other applications of this monodromy-based approach, including determining identifiability in biological models \cite{BHM}.

\section*{Acknowledgements}
The first three authors thank the Simons Institute for the Theory of Computing in Berkeley, CA for their generous support while in residence during the program on {\em Algorithms and Complexity in Algebraic Geometry}.
JDH was also supported by DARPA YFA, NSF DMS-1262428, NSF ACI-1460032,
and a Sloan Research Fellowship. 
GO is member of GNSAGA.
AJS was partially supported by NSF ACI-1440607.

\end{document}